\documentclass{article}

\usepackage{amssymb, amsmath}
  \newtheorem{lemma}{Lemma}%

\newtheorem{remark}{Remark}

\newcommand{\vep}{\varepsilon}

\newcommand{\QED}{\hfill
$\underline{\underline{QED}}$}

\newcommand{\la}{\lambda}

\newcommand{\beq}{\begin{equation}}
\newcommand{\eeq}{\end{equation}}
\newcommand{\beqn}{\begin{eqnarray}}
\newcommand{\eeqn}{\end{eqnarray}}
\newcommand{\beqnn}{\begin{eqnarray*}}
\newcommand{\eeqnn}{\end{eqnarray*}}
\newcommand{\nn}{|\!|}

\newcommand{\rr}{\mathbb{R}}
\newcommand{\n}{\mathbb{N}}

\begin{document}
\title{Existence result for a Neumann problem}

\author{Nikolaos Halidias \\
University of the Aegean \\ Department of Statistics and Actuarial
Science
\\ Karlovassi, 83200 \\ Samos \\ Greece \\ email: nick@aegean.gr}

\maketitle

%\subjclass{2000 Mathematics Subject Classification. 32J15,34J89}%

\begin{abstract} In this paper
we are going to show the existence of a nontrivial solution to the
following model problem,
\begin{equation*}
\left\{
\begin{array}{lll}
-\Delta (u)  = 2uln(1+u^2)+\frac{|u|^2}{1+u^2}2u+usin(u) \mbox{
a.e. on } \Omega
\\ \frac{\partial u}{\partial \eta} = 0
 \mbox{ a.e. on } \partial \Omega.
\end{array}
\right.
\end{equation*}
As one can see the right hand side is superlinear. But we can not
use an Ambrosetti-Rabinowitz condition in order to obtain that the
corresponding energy functional satisfies (PS) condition. However,
it follows that the energy functional satisfies the Cerami (PS)
condition.\footnote{2000 Mathematics Subject Calssification:
35A15, 35J20, 35J25
 \\Keywords: Mountain-Pass Theorem, critical point, Cerami $PS$ condition.}
\end{abstract}

\section{Introduction}In this paper
we are going to show the existence of a nontrivial solution to the
following model problem,
\begin{equation}
\left\{
\begin{array}{lll}
-\Delta (u)  = 2uln(1+u^2)+\frac{|u|^2}{1+u^2}2u+usin(u) \mbox{
a.e. on } \Omega
\\ \frac{\partial u}{\partial \eta} = 0
 \mbox{ a.e. on } \partial \Omega.
\end{array}
\right.
\end{equation}
As one can see the right hand side is superlinear. But we can not
use an Ambrosetti-Rabinowitz condition in order to obtain that the
corresponding energy functional satisfies (PS) condition. Let us
recall the well-known Ambrosetti-Rabinowitz condition:

 There
exists some $\theta > 2$ such that
\begin{eqnarray*}
0 < \theta F(u) \leq F^{'}(u)u,
\end{eqnarray*}
for all $|u| > M$ for big enough $M$.

We can see that for our model problem there is not such a $\theta
> 2$. But, we will show that there exists a sequence $\theta_n > 2$
with $\theta_n \to 2$ which have the desired effectiveness for our
problem when we use the Cerami (PS) condition.

For such kind of problems there are some papers that extends the
well-known Ambrosetti-Rabinowitz condition. For example one can
also see the very interesting result of D.G. de Figueiredo-J. Yang
\cite{Fig} who considers semilinear problems such that the
corresponding energy functional does not satisfy a (PS) condition.
Also, Gongbao Li-HuanSong Zhou \cite{LIZHOU} made some progress in
this direction. But they assume that $f(u) \geq 0$ for all $u \in
\rr$. Take problem (1) and see that  that $f (u) \to -\infty$ as
$u \to -\infty$, thus we can not say that $f (u) \geq 0$ for all
$u \in \rr$. So, we can not use their method in order to obtain a
nontrivial solution. Finally, let us mention the work of
Costa-Magalhaes \cite{CMAG}. In this paper we extend the results
of \cite{CMAG} using deferent arguments in our proof. The authors
there have proposed the following hypotheses, among others,
\begin{eqnarray*}
f(s)s-p F(s) \geq a |s|^{\mu}, \mbox{ for all } s \in \rr,
\end{eqnarray*}
with $\mu \geq \frac{N/p}{q-p}$ for some $q < p^* =
\frac{Np}{N-p}$. Thus $f(s)s-pF(s)$ must grows faster than
$s^\frac{N/p}{q-p}$. Here, we extend this result because we do not
need such a bound for $\mu$.

Our existence theorem, considers more general Neumann problems
than our model problem and at the end of the paper we give a
second example. We must also note that non of the above papers
considers Neumann problems.

Let us mention some facts that we are going to use later. It is
well known that $W^{1,p}(\Omega) = \rr \oplus W$ with $W = \{ u
\in W^{1,p}(\Omega): \int_{\Omega} u(x) dx = 0 \}$. We can
introduce the following number,
\begin{eqnarray*}
\la_1 = \inf \{ \frac{\nn Dw \nn_p^p}{\nn w \nn^p_p}:w \in W, w
\neq 0 \}.
\end{eqnarray*}
From Papalini \cite{Papalini} we know that $\la_1 >0$ and if $w
\in W$ is such that $\nn w \nn_p = 1$, $\nn Dw \nn_p = \la_1$ then
$w$ is an eigenfunction of the following problem,
\begin{equation}
\left\{
\begin{array}{lll}
-\Delta_p(u)  = \la_1 |u|^{p-2}u \mbox{ a.e. on } \Omega
\\ \frac{\partial u}{\partial \eta_p} = 0
 \mbox{ a.e. on } \partial \Omega, \;\; 2 \leq p < \infty.
\end{array}
\right.
\end{equation}

Let us introduce the $(PS)$ that we are going to use. \\
{\bf Cerami $(PS)$ condition} Let $X$ be a Banach space and $I:X
\to \rr$.
 For every $\{ u_n \} \subseteq W^{1,p}(\Omega)$ with $|I(u_n)|
 \leq M$ and $(1+\nn u_n \nn_{1,p})<I^{'}(u_n),\phi> \to 0$ for
 every $\phi \in W^{1,p}(\Omega)$ there exists a strongly
 convergent subsequent. This condition has introduced by Cerami
 (see \cite{Cerami}, \cite{BBF}).

\section{Basic Results}
We are going to show an existence result for the following Neumann
problem,
\begin{equation}
\left\{
\begin{array}{lll}
-\Delta_p (u)  = f(x,u) \mbox{ a.e. on } \Omega
\\ \frac{\partial u}{\partial \eta_p} = 0
 \mbox{ a.e. on } \partial \Omega, \; p \geq 2.
\end{array}
\right.
\end{equation}
We suppose that $\Omega$ is a bounded domain with sufficient
smooth boundary $\partial \Omega$. By $\Delta_p$ we denote the
well-known $p$-Laplacian operator, i.e. $\Delta_p (u) = div(\nn Du
\nn^{p-2}Du)$.

From now on we will denote by $F(x,u) = \int_o^r f(x,r)dr$ and
$h(x,u) = \frac{F(x,u)}{|u|^p}$. We suppose  the following
assumptions on $f$,

{\bf H(f)}$f: \Omega \times \rr \to \rr$ is a Carath\'eodory
function such that
\begin{enumerate}
\item[(i)] for almost all $x \in \Omega$ and for all $u \in \rr$ we have that $|f(x,u)| \leq
C(1+|u|^{\tau-1})$, with $\tau < p^* = \frac{np}{n-p}$ and $C>0$;

\item[(ii)] uniformly for almost all $x \in \Omega$
$$\liminf_{|r| \to \infty} h(x,r) \geq \mu > 0;$$

\item[(iii)]  there exists some $r < p$ such that for every $k \in
L^{p-r}(\Omega)$,we have
\begin{eqnarray*}
\lim_{|u| \to \infty} \frac{p F(x,u) - f(x,u) u}{k(x)
|u|^{p-r}h(x,u)} = -\infty.
\end{eqnarray*}
 for almost all $x \in \Omega$;
\item[(iv)] uniformly for almost all $ x \in \Omega$ we have
 $\limsup_{u \to 0} h(x,u) \leq \theta (x)$, with $\theta (x) \leq
 \frac{\la_1}{p}$ and $\int_{\Omega} (\la_1-\theta (x))|w(x)|^pdx
 > 0$ for every $w$  an eigenfunction corresponding to the
 second eigenvalue $\la_1$.

\end{enumerate}

\begin{remark} Condition $H(f)(iii)$ means that there exists some
$r < p$ such that for every $k \in L^{p-r}(\Omega)$ and every
$C>0$ we can find big enough $M>0$ such that
\begin{eqnarray*}
C k(x)|u|^{p-r}h(x,u) + p F(x,u) - f(x,u) u \leq 0.
\end{eqnarray*}
\end{remark}

Let us define first the energy functional $I:W^{1,p}(\Omega) \to
\rr$ by $I(u) = \frac{1}{p} \nn Du \nn_p^p - \int_{\Omega}
F(x,u(x))dx$. Under  conditions $H(f)$ it is well known that $I$
is well defined and a $C^1$ functional. We are going to use the
Mountain-Pass Theorem, so our first lemma
 is that $I$ satisfies the Cerami $(PS)$ condition.

\begin{lemma}
$I$ satisfies the Cerami $(PS)$ condition.
\end{lemma}

{\bf Proof}

Let $\{ u_n \} \subseteq W^{1,p}(\Omega)$ such that $|I(u_n)| \leq
M$ and $(1+\nn u_n \nn_{1,p})<I^{'}(u_n),\phi> \to 0$ for every
$\phi \in W^{1,p}(\Omega)$. We must show that $u_n$ is bounded.
Suppose that $\nn u_n \nn_{1,p} \to \infty$. We will show that
$\nn D u_n \nn_p \to \infty$. Indeed, from the choice of the
sequence we have

\begin{eqnarray*}
-M \leq \nn D u_n \nn_p^p - \int_{\Omega} p F(x,u_n)dx
\leq M, \Rightarrow \\
\int_{\Omega} F(x,u_n)dx \leq M + \nn D u_n \nn_p^p \Rightarrow \\
c \nn u_n \nn_p^p \leq M + \nn D u_n \nn_p^p,
\end{eqnarray*}
here have used $H(f)(ii)$.

 Thus, we can not suppose that $\nn D u_n
\nn_p$ is bounded because then $\nn u_n \nn_p \to \infty$ and then
from the above relation we obtain a contradiction. So, it follows
that $\nn D u_n \nn_p \to \infty$ and moreover
\begin{eqnarray*}
\nn u_n \nn_p^p \leq (\vep_n + c)\nn D u_n \nn_p^p,
\end{eqnarray*}
with $\vep_n \to 0$.

It follows then that  there exists some $c_1,c_2$ such that

\begin{eqnarray}
c_1 \nn u_n \nn_{1,p} \leq \nn D u_n \nn_p \leq c_2 \nn u_n
\nn_{1,p}.
\end{eqnarray}

Then, from the choice of the sequence it follows
\begin{eqnarray}
-M \leq -\nn Du_n \nn^p_p + \int_{\Omega} p F(x,u_n) dx \leq M,
\end{eqnarray}
and choosing $\phi = u_n$
\begin{eqnarray}
-\vep_n \frac{\nn u_n \nn_{1,p}}{1+\nn u_n \nn_{1,p}} \leq \nn
Du_n \nn_p^p - \int_{\Omega} f(x,u_n) u_ndx \leq \vep_n \frac{\nn
u_n \nn_{1,p}}{1+\nn u_n \nn_{1,p}}.
\end{eqnarray}
Consider now the sequence  $a_n = \frac{1}{p \nn u_n
\nn_{1,p}^r}$. Then multiply inequality (5) with $a_n+1$,
substituting with (6) and using (4) we arrive at
\begin{eqnarray}
& & c \nn u_n \nn_{1,p}^{p-r} \leq  a_n \nn Du_n \nn_p^p \leq \nonumber \\
& & \int_{\Omega}(a_n+1)p F(x,u_n) - f(x,u_n)u_ndx + (a_n+1)M
+\vep_n \frac{\nn u_n \nn_{1,p}}{1+\nn u_n \nn_{1,p}}.
\end{eqnarray}

Let $y_n(x) = \frac{u_n(x)}{\nn u_n \nn_{1,p}}$. Then, it is clear
that there exists some $k \in L^p(\Omega)$ such that $|y_n(x)|
\leq k(x)$ a.e. on $\Omega$.

Let $\Omega_1 = \{ x \in \Omega:|u(x)| \leq M \}$. In view of
$H(f)(iv)$ we have that
\begin{eqnarray*}
\int_{\Omega_1} (a_n+1)p F(x,u_n(x)) - f(x,u_n(x))u_n(x)dx\leq C,
\end{eqnarray*}
for every $M>0$. Choosing big enough $M>0$ we can estimate,
\begin{eqnarray*}
& & \int_{\Omega \setminus \Omega_1} (a_n+1)p F(x,u_n(x)) - f(x,u_n(x))u_n(x)dx\leq \\
& & \int_{\Omega \setminus \Omega_1}  |y_n(x)|^r |u_n(x)|^{p-r}
h(x,u_n(x))+p F(x,u_n(x)) - f(x,u_n(x))u_n(x)dx  \\ & & \leq
\int_{\Omega \setminus \Omega_1}k(x)|u_n(x)|^{p-r}h(x,u_n(x)) + p
F(x,u_n(x))-f(x,u_n(x))u_n(x) dx \leq 0.
\end{eqnarray*}

Going back to (7), we obtain a contradiction to the hypothesis
that $u_n$ is not bounded.

Finally, using well-known arguments we can prove that in fact $\{
u_n \}$ have a convergent subsequence.

\QED

\begin{lemma} There exists some $\xi \in \rr$  such
that $I(\xi ) \leq 0$.
\end{lemma}

{\bf Proof}

We claim that there exists big enough $\xi \in \rr$ such that
$I(\xi) \leq 0$.  Suppose not. Then there exists a sequence $\xi_n
\to \infty$ such that $I(\xi_n ) \geq c > 0$. That means
\begin{eqnarray*}
- \int_{\Omega} F(x,\xi_n  )dx \geq c > 0.
\end{eqnarray*}
Using now $H(f)(ii)$ we can say that for almost all $x \in \Omega$
and all $u \in \rr$ we have that $F(x,u) \geq \mu |u|^p -c$. So,
it follows that
\begin{eqnarray*}
\mu |\xi_n|^p  \leq c \mbox{ for every } n \in \n.
\end{eqnarray*}
But this is a contradiction.

\QED

\begin{lemma} There exists some $\rho >0$ small enough and $a >0$ such that
$I(u) \geq a$ for all $\nn u \nn_{1,p} = \rho$ with $u \in W$.
\end{lemma}

{\bf Proof} Suppose that this is not true. Then there exists a
sequence $\{ u_n \} \subseteq W$ such as $\nn u_n \nn_{1,p} =
\rho_n$ with $\rho_n \to 0$, with the property that $I(u_n) \leq
0$. So we arrive at
\begin{eqnarray}
\nn Du_n \nn_p^p \leq p\int_{\Omega} F(x,u_n(x))dx
\end{eqnarray}
 Let $y_n(x) = \frac{u_n(x)}{\nn u_n \nn_{1,p}}$. Using $H(f)(i),(iv)$ we
can prove that there exists $\gamma > 0$ such that
\begin{eqnarray*}
p F(x,u) \leq (\theta (x) + \vep)|u|^p + \gamma |u|^{p^*}
\end{eqnarray*}
Take in account the last estimation and dividing (8) with $\nn u_n
\nn_{1,p}^p$ we arrive at
\begin{eqnarray}
\la_1 \nn y_n \nn_p^p \leq \nn Dy_n \nn_p^p  \leq \int_{\Omega}
(\theta (x) + \vep)|y_n(x)|^pdx + \gamma_1 \nn u_n \nn_p^{p^*-p}.
\end{eqnarray}
Recall that $y_n \to y$ strongly in $L^p(\Omega)$. Using the lower
semicontinuity of the norm we arrive at $\nn Dy \nn_p \leq \la_1
\nn y \nn_p$  and from the variational characterization of the
second eigenvalue these quantities are in fact equal. Note that
$y_n \to y$ weakly in $X$ and recall that $\nn Dy_n \nn_p \to
\la_1 \nn y \nn_p = \nn Dy \nn_p$. Then from the uniform convexity
of $X$ we have $y_n \to y$ strongly  in $X$ and $y \neq 0$. Thus
$y \in W$ is an eigenfunction of $(-\Delta_p, W)$.

Going back to (9) and taking the limit we arrive at
\begin{eqnarray*}
\int_{\Omega} (\la_1 - \theta (x))|y(x)|^pdx \leq 0.
\end{eqnarray*}
But this is a contradiction.

 \QED

Then the existence of a nontrivial solution for problem (3)
follows from a variant of  Mountain-Pass Theorem (see Struwe
\cite{Struwe}, Thm. 8.4 and Example 8.2, or \cite{BBF}, Thm. 2.3
and Prop. 2.1).

\section{Applications to Differential Equations}

Consider the following elliptic equation,
\begin{equation}
\left\{
\begin{array}{lll}
-\Delta (u)  = 2uln(1+u^2)+\frac{|u|^2}{1+u^2}2u+usin(u) \mbox{
a.e. on } \Omega
\\ \frac{\partial u}{\partial \eta} = 0
 \mbox{ a.e. on } \partial \Omega.
\end{array}
\right.
\end{equation}
Here, as before, $\Omega \subseteq \rr^n$ is a bounded domain with
smooth enough boundary $\partial \Omega$.

We can check that the corresponding energy functional does not
satisfy an Ambrosetti-Rabinowitz type condition. Moreover, we can
not say that $f(u) \geq 0$ nor that $f(u)+f(-u) = 0$, thus we can
not use the arguments of \cite{LIZHOU}, \cite{LIZHOU2}, even if
our problem had Dirichlet boundary conditions. Finally, $f(\cdot)$
does not satisfy the condition of \cite{Fig} because does not
exist $p>1$ such that $f(u) \geq \mu u^p$ for all $u \geq T$ for
big enough $T$. Also, we can choose big enough $n \in \n$ (i.e.
the dimension of our problem) and see that the above problem does
not satisfy the conditions of \cite{CMAG}.

However, we can check that $f(\cdot)$ satisfies the conditions
that we have proposed.

We can see also that $h$ did not have to go to infinity. Take for
example as $F(u) = |u|^p (a+(b-a)\frac{|u|^r}{1+|u|^r})$. Choose
$a < \frac{\la_1}{p}, b > 0$ and for a suitable choice of $r$ (for
example $r < p$) we can see that $f$ satisfies all the hypotheses
that we have proposed.

\end{document}